\documentclass[12pt]{article}

\usepackage{fullpage}

\usepackage[utf8]{inputenc}

\usepackage[backend=biber, style=numeric-comp, giveninits=true]{biblatex}
\usepackage{listings}
\usepackage{varioref}
\usepackage{tkz-berge}

\lstdefinestyle{mystyle}{
    backgroundcolor=\color{gray!15},   
    commentstyle=\color{black},
    keywordstyle=\color{black},
    numberstyle=\tiny\color{black},
    stringstyle=\color{black},
    basicstyle=\ttfamily\footnotesize,
    breakatwhitespace=false,         
    breaklines=true,                 
    captionpos=b,                    
    keepspaces=true,                 
    numbers=left,                    
    numbersep=5pt,                  
    showspaces=false,                
    showstringspaces=false,
    showtabs=false,                  
    tabsize=2
}

\usepackage{xcolor} 
\usepackage{url}
\usepackage{enumerate} 

\usepackage{mathrsfs}

\usepackage{amsmath,amssymb,amsthm} 

\newtheorem{theorem}{Theorem}
\newtheorem{proposition}{Proposition}
\newtheorem{lemma}{Lemma}
\newtheorem{corollary}{Corollary}
\theoremstyle{remark}
\newtheorem{remark}{Remark}
\newtheorem{example}{Example}
\theoremstyle{definition}
\newtheorem{definition}{Definition}

\DeclareMathOperator{\Aut}{Aut}
\DeclareMathOperator{\Rot}{Rot}
\DeclareMathOperator{\Kappa}{K}
\DeclareMathOperator{\Fix}{Fix}
\DeclareMathOperator{\Wheel}{W}

\newcommand{\stern}[1]{\Kappa_{1,#1}}
\newcommand{\AG}[1]{\Aut(\Gamma_{#1})}
\newcommand{\RG}[1]{\Rot(\Gamma_{#1})}

\author{Klaus Dohmen and Mandy Lange-Geisler\\ Department of Mathematics\\
  Mittweida University of Applied Sciences\\ 09648 Mittweida, Germany}

\title{\bf The Orbital Bivariate Chromatic Polynomial}

\begin{document}

\maketitle

\begin{abstract}
\emph{Abstract.}  We introduce the \emph{orbital bivariate chromatic
  polynomial}, which counts generalized vertex colorings of a graph up to the
action of a group of symmetries. More precisely, adjacent vertices must either
receive distinct colors from a set of $\lambda$ colors or may receive the same
color from a distinguished subset of $\lambda-\mu$ colors.  This new graph
polynomial simultaneously generalizes the orbital chromatic polynomial due to
Cameron and Kayibi (2007) and bivariate chromatic polynomial due to Dohmen,
P\"onitz, and Tittmann (2003).  We establish fundamental structural properties
and derive explicit formulas for several graph families, including complete
graphs, complete bipartite graphs, paths, cycles, and wheels.  In particular,
several of these results are new even in the special case of the orbital
chromatic polynomial.  For cycles, our results generalize Moreau's necklace
polynomial and lead to a new generalization and proof of Fermat's Little
Theorem.
\par\medskip

\par\medskip

\end{abstract}

\section{Introduction}

Since its introduction by Birkhoff in 1912 in connection with the four-color
problem \cite{Birkhoff:1912}, the chromatic polynomial has become a subject of
considerable interest in combinatorial mathematics.  For a finite, undirected
graph, its chromatic polynomial counts the number of vertex colorings with a
prescribed number of colors in which adjacent vertices receive distinct
colors. For an accessible introduction to this subject,
the reader is referred to Read \cite{Read:1968}.

Among the many variants and generalizations of the chromatic polynomial, the
\emph{bivariate chromatic polynomial}, introduced by Dohmen, P\"onitz, and
Tittmann \cite{Dohmen:2003}, has received considerable attention
\cite{Averbouch:2010, Beck:2023, Trinks:2012, Gerling:2019}. This two-variable
polynomial generalizes not only the chromatic polynomial but also the
independence polynomial, and it is closely related to the matching
polynomial. By definition, the bivariate chromatic polynomial
$P_{\Gamma}(\lambda,\mu)$ of a finite, undirected graph $\Gamma$ counts the
colorings of the vertices of $\Gamma$ with $\lambda$ colors such that adjacent
vertices receive distinct colors from a set of $\lambda$ colors, or the same
color from a distinguished subset of $\lambda-\mu$ colors (i.e., colors larger
than $\mu$).

In some applications---such as chemistry---colorings that are related by a
graph symmetry should be regarded as equivalent.  To address this, Cameron and
Kayibi introduced the \emph{orbital chromatic polynomial} in
\cite{Cameron:2007}, which counts proper $\lambda$-colorings of a finite,
undirected graph up to the action of a subgroup of its automorphism group.  In
\cite{Cameron:2007, Cameron:2008, Kim:2014}, the orbital chromatic polynomial
and its roots have been studied for small graphs, as well as for the classes
of edgeless and complete graphs.

While the orbital chromatic polynomial extends the chromatic polynomial to
graphs with symmetries, its computation via Burnside's lemma requires quotient
graphs, which may contain loops. Since a loop prevents proper colorings, the
corresponding terms vanish, and structural information encoded in such
quotient graphs is lost. The bivariate chromatic polynomial avoids this by
introducing an additional parameter $\mu$, which permits certain otherwise
improper colorings to be counted. Consequently, quotient graphs with loops may
contribute nontrivially, thereby retaining structural information that is not
captured by the orbital chromatic polynomial. This motivates the introduction
of a bivariate analogue.

Our new \emph{orbital bivariate chromatic polynomial}, presented in this
article, simultaneously generalizes the orbital chromatic
polynomial~\cite{Cameron:2007} and the bivariate chromatic
polynomial~\cite{Dohmen:2003}.  It counts the same colorings as the bivariate
chromatic polynomial, up to the same equivalence relation induced by graph
symmetries as in the orbital chromatic polynomial.  After establishing
fundamental properties of this new invariant, we derive general formulas for
its evaluation on several families of graphs, including complete graphs,
complete bipartite graphs, paths, cycles, and wheels. In particular, some of
these expansions are new even for the orbital chromatic polynomial. As a
by-product, we obtain a new generalization and proof of Fermat's Little
Theorem.

Our article is organized as follows. In Section \ref{sec:bivariate} we recall
the definition of the bivariate chromatic polynomial introduced
in~\cite{Dohmen:2003}, together with closed-form formulas, including new
formulas involving Fibonacci and Lucas sequences.  In
Section~\ref{sec:orbital} we introduce the orbital bivariate chromatic
  polynomial, which simultaneously generalizes the orbital chromatic
polynomial \cite{Cameron:2007} and the bivariate chromatic
polynomial~\cite{Dohmen:2003}, and establish some of its fundamental properties.
In Section~\ref{sec:spec-graph-class} we derive general formulas for edgeless
graphs, complete graphs, complete bipartite graphs, stars, paths, cycles, and
wheels of arbitrary order.  In Section~\ref{sec:final-treat}, we relate the
orbital bivariate chromatic polynomial to Moreau's necklace polynomial, which
arises as a specialization of our polynomial.  Finally, using our results for
cycles, we obtain a new generalization and proof of Fermat's Little Theorem.

Throughout this article, all graphs are assumed to be finite and undirected,
and they may contain loops and multiple edges.  Since we consider both graphs
and groups, we denote graphs by capital Greek letters and groups by capital
Roman letters.

\section{The bivariate chromatic polynomial}
\label{sec:bivariate}

In this section, we recall the bivariate chromatic polynomial introduced in
\cite{Dohmen:2003}. Besides its definition and basic properties, we present an
illustrative example and several known formulas for graph families that will
play a role in later sections.

\subsection{Notions and basic definition}

For every graph $\Gamma$ and each $\lambda\in\mathbb{N}$, a
\emph{$\lambda$-coloring} of $\Gamma$ is a mapping $f$ from the vertex set of
$\Gamma$ to $\{1,\dots,\lambda\}$.  For $\mu=0,\dots,\lambda$ we call $f$
\emph{$\mu$-proper} if, for every pair of adjacent vertices $v$ and $w$ of
$\Gamma$, either $f(v)\neq f(w)$, or $f(v)=f(w)>\mu$; that is, adjacent
vertices receive different colors from $\{1,\dots,\lambda\}$, or the same
color from $\{\mu+1,\dots,\lambda\}$.

The classical \emph{chromatic polynomial} $P_{\Gamma}(\lambda)$, introduced by
Birkhoff \cite{Birkhoff:1912}, counts the $\mu$-proper $\lambda$-colorings
of~$\Gamma$ with $\mu=\lambda$. The following definition generalizes this
concept.

\begin{definition}[\cite{Dohmen:2003}]
For every graph $\Gamma$, every $\lambda\in\mathbb{N}$ and
$\mu=0,\ldots,\lambda$, we let $P_{\Gamma}(\lambda,\mu)$ denote the number of
$\mu$-proper $\lambda$-colorings of $\Gamma$.  We refer to
$P_{\Gamma}(\lambda,\mu)$ as the \emph{bivariate chromatic polynomial} of~$G$.
\end{definition}

It follows from the results in \cite{Dohmen:2003} that
$P_{\Gamma}(\lambda,\mu)$ is indeed a polynomial in $\lambda$ and~$\mu$.  For
$\lambda=\mu$ it coincides with the chromatic polynomial $P_\Gamma(\lambda)$,
while for $\mu=0$ it coincides with $\lambda^{n(\Gamma)}$, where $n(\Gamma)$
is the number of vertices of~$\Gamma$.

The bivariate chromatic polynomial has been determined for several families of
graphs, e.g., complete graphs, complete bipartite graphs, paths, and
cycles~\cite{Dohmen:2003}.  As shown in~\cite{Dohmen:2003}, it can be computed
in polynomial time for graphs of bounded treewidth.  Further graph classes
have been studied in~\cite{Gerling:2019}.

\subsection{An illustrative example}
\label{sec:an-illustr-example}

The following example illustrates the combinatorial interpretation of the
bivariate chromatic polynomial, with particular emphasis on the role of~$\mu$.

\newcommand{\exgra}{\Gamma_4+e}

\begin{example}

\label{example:bivar-chrom-polyn-1}
\begin{figure}[!t]
\centering
\begin{tikzpicture}[scale=.4]
\SetGraphUnit{4}
\GraphInit[vstyle=Classic]
\begin{scope}[rotate=45]
\Vertices[Lpos=45]{circle}{2,1,4,3}
\end{scope}
\Edges(4,1,2,3,4,2)
\end{tikzpicture}
\caption{Example graph}
\label{fig:example}
\end{figure}

Consider the graph $\exgra$, obtained by adding a diagonal to a $4$-cycle, as
depicted in Figure~\vref{fig:example}.
\begin{table}[ht]
\centering
\begin{tabular}{p{0.45\linewidth} c c}
\hline
monochromatic edges & choices & $\mu$-proper $\lambda$-colorings \\
\hline
none
& $1$
& $\lambda(\lambda-1)(\lambda-2)^2$ \\
only the diagonal
& $1$
& $(\lambda-\mu)(\lambda-1)^2$ \\
exactly one edge of the $4$-cycle
& $4$
& $4(\lambda-\mu)(\lambda-1)(\lambda-2)$ \\
exactly two opposite edges of the $4$-cycle
& $2$
& $2(\lambda-\mu)(\lambda-\mu-1)$ \\
exactly two adjacent edges of the $4$-cycle
& $4$
& $4(\lambda-\mu)(\lambda-1)$ \\
all edges
& $1$
& $\lambda-\mu$ \\
\hline
\end{tabular}
\caption{Number of $\mu$-proper $\lambda$-colorings of $\exgra$
w.r.t.\ monochromatic edges}
\label{tab:bivariate-example}
\end{table}
In Table~\vref{tab:bivariate-example}, we determine its bivariate chromatic
polynomial by counting $\mu$-proper $\lambda$-colorings according to which
edges are monochromatic.  Summing the entries in the third row gives
\begin{gather*}
P_{\exgra}(\lambda,\mu) = \lambda^4 - 5\lambda^2\mu +6\lambda\mu + 2\mu^2 -4\mu.
\end{gather*}
For $\lambda=\mu$, this specializes to the chromatic polynomial
$\lambda(\lambda-1)(\lambda-2)^2$ of~$\exgra$.
\end{example}

\subsection{Explicit formulas}
\label{sec:explicit-formulas}

The following propositions give explicit expressions for the bivariate
chromatic polynomial of complete graphs, complete bipartite graphs, stars,
paths, and cycles. In Section~\ref{sec:spec-graph-class}, we use these
expressions to derive the corresponding formulas for the orbital bivariate
chromatic polynomial of these graph classes.

In the sequel, we use $\Kappa_n$ to denote the complete graph on $n$ vertices,
$\Kappa_{m,n}$ to denote the complete bipartite graph on $m+n$ vertices,
$\Pi_n$ to denote the path of length $n$ (having $n+1$ vertices), and
$\Gamma_n$ to denote the cycle of length $n$.  We will further clarify these
notations in Section~\ref{sec:spec-graph-class}, when referring to the
vertex and edge sets of these graphs.

We use $x^{\underline{k}}$ to denote the $k$-th falling factorial of~$x$;
i.e., $x^{\underline{k}} = x(x-1)\dots(x-k+1)$, and
$\left\{ k \atop l \right\}$ to denote the Stirling numbers of the second
kind.

\begin{proposition}[\cite{Dohmen:2003}] For every $n,\lambda\in\mathbb{N}$ and
each $\mu=0,\dots,\lambda$, we have
\begin{gather}
\label{eq:8}
P_{\Kappa_n}(\lambda,\mu) = \sum_{k=0}^n {n\choose k} (\lambda-\mu)^k
\mu^{\underline{n-k}}.
\end{gather}
\end{proposition}

\begin{proposition}[\cite{Dohmen:2003}] For every $m,n,\lambda\in\mathbb{N}$
and each $\mu=0,\dots,\lambda$, we have
\begin{gather}
\label{eq:19}
P_{\Kappa_{m,n}}(\lambda,\mu) = \sum_{k=0}^m {m \choose k} (\lambda-\mu)^{m-k}
\sum_{l=0}^k \left\{ k\atop l\right \}  (\lambda-l)^n \mu^{\underline{l}} \,  .
\end{gather}
In particular,
\begin{equation}
\label{eq:3}
P_{\stern{n}}(\lambda,\mu) = \lambda^n(\lambda-\mu) + (\lambda-1)^n\mu .
\end{equation}
\end{proposition}

After the case of stars $K_{1,n}$, the bivariate chromatic polynomial of paths
and cycles can also be expressed in closed form.

\begin{proposition}[\cite{Dohmen:2015}] For every $n\in\mathbb{N}_0$, every
$\lambda\in\mathbb{N}$ and each $\mu=0,\dots,\lambda$, except for
$\lambda=\mu=1$, we have
\label{prop:path}
\begin{equation}
\label{eq:6}
P_{\Pi_n}(\lambda,\mu) = \frac{\sqrt{D}-\lambda-1}{2\sqrt{D}}\left(
\frac{\lambda-1-\sqrt{D}}{2}\right)^{n+1} \! + \,
\frac{\sqrt{D}+\lambda+1}{2\sqrt{D}}\left( \frac{\lambda-1+\sqrt{D}}{2} \right)^{n+1},
\end{equation}
where $D=(\lambda+1)^2 - 4\mu$. Moreover, $P_{\Pi_0}(1,1)=1$ and
$P_{\Pi_n}(1,1)=0$ for $n\ge 1$.
\end{proposition}

\begin{proposition}[\cite{Dohmen:2015}] For every $n,\lambda\in\mathbb{N}$ and
each $\mu=0,\dots,\lambda$, we have 
\label{prop:cycle}
\begin{equation}
\label{eq:5}
P_{\Gamma_n}(\lambda,\mu) = \left(
\frac{\lambda-1-\sqrt{D}}{2}\right)^{n} \! + \,
\left( \frac{\lambda-1+\sqrt{D}}{2} \right)^{n} + (-1)^n(\mu-1)\, ,
\end{equation}
where $D=(\lambda+1)^2 - 4\mu$.
\end{proposition}

\subsection{Connection with Lucas sequences}
\label{sec:conn-with-lucas}

The following connection of Propositions \ref{prop:path} and \ref{prop:cycle}
with the Lucas sequences appears to be new and will be used later on in
Subsection~\ref{sec:rg}. Recall from \cite{Lucas:1878} that for
$P,Q\in \mathbb{Z}$ and $n>1$, the Lucas sequences of the first and second
kind satisfy the recurrence relations
\begin{align*}
  U_n(P,Q) & = P\cdot U_{n-1}(P,Q) - Q\cdot U_{n-2}(P,Q), \\
  V_n(P,Q) & = P\cdot V_{n-1}(P,Q) - Q\cdot V_{n-2}(P,Q),
\end{align*}             
with initial conditions $U_0(P,Q) = 0$, $U_1(P,Q) = 1$, $V_0(P,Q) = 2$,
$V_1(P,Q) = P$.  It is well-known and straightforward to check that if
$P^2 > 4Q$, then
\begin{align}
  U_n(P,Q) & = \frac{1}{\sqrt{D}}\left(\left(\frac{P+\sqrt{D}}{2}\right)^n -
             \left(\frac{P-\sqrt{D}}{2}\right)^n \right), \notag \\
  V_n(P,Q) & = \left(\frac{P+\sqrt{D}}{2}\right)^n +
             \left(\frac{P-\sqrt{D}}{2}\right)^n,
             \label{eq:25}
\end{align}
where $D = P^2-4Q$ (cf.~\cite{Dickson:2005}).  Thus, with $P=\lambda-1$
and $Q = \mu-\lambda$, Eqs.~(\ref{eq:6}) and (\ref{eq:5}) in
Propositions \ref{prop:path} and \ref{prop:cycle} can equivalently be stated
as
\begin{align}
  P_{\Pi_n}(\lambda,\mu) & = U_{n+2}(\lambda-1,\mu-\lambda) +
                           U_{n+1}(\lambda-1,\mu-\lambda) \quad (\text{except for $\lambda=\mu=1$}),
                           \tag{\ref{eq:6}'} \label{eq:61} \\
  P_{\Gamma_n}(\lambda,\mu) & = V_n(\lambda-1, \mu-\lambda) + (-1)^n(\mu-1)\,
                              \tag{\ref{eq:5}'}. \label{eq:51}
\end{align}
There is also a connection with Fibonacci and Lucas polynomials,
as defined by
\begin{align*}
  F_n(x) & = x F_{n-1}(x) + F_{n-2}(x) \quad (n>1); \quad F_0(x)=0, \quad
           F_1(x)=1; \\
  L_n(x) & = x L_{n-1}(x) + L_{n-2}(x) \quad (n>1); \quad L_0(x)=2, \quad L_1(x)=x.
\end{align*}
Since $F_n(x) = U_n(x,-1)$ and $L_n(x) = V_n(x,-1)$, we conclude from
Eqs.~\eqref{eq:61} and \eqref{eq:51} that
\begin{align*}
  P_{\Pi_n}(\lambda,\lambda-1) & = F_{n+2}(\lambda-1) + F_{n+1}(\lambda-1), \\
   P_{\Gamma_n}(\lambda,\lambda-1) & = L_n(\lambda-1) + (-1)^n(\lambda-2) .
\end{align*}
In particular, $P_{\Pi_n}(2,1)$ gives the $(n+3)$-rd Fibonacci number
$F_{n+3} = F_{n+3}(1)$, while $P_{\Gamma_n}(2,1)$ gives the $n$-th Lucas number
$L_n = L_n(1)$.

\section{The orbital bivariate chromatic polynomial}
\label{sec:orbital}

In this section, we introduce the orbital bivariate chromatic polynomial,
derive a formula for it via Burnside’s lemma, and establish several of its
basic properties.

\subsection{Definition and Burnside formula}

With each graph $\Gamma$ and each subgroup $G$ of its automorphism group
$\Aut(\Gamma)$, we associate a two-variable function in $\lambda$ and $\mu$
that counts the non-equivalent $\mu$-proper $\lambda$-colorings of $\Gamma$
under the action of~$G$. This function is indeed a polynomial in $\lambda$ and
$\mu$, which forms the subject of our investigation.

\begin{definition}
\label{def:orbital}
For every graph $\Gamma$, every subgroup $G$ of $\Aut(\Gamma)$, every
$\lambda\in\mathbb{N}$ and $\mu=0,\ldots,\lambda$, we define
$OP_{\Gamma,G}(\lambda,\mu)$ as the number of equivalence classes of
$\mu$-proper $\lambda$-colorings of $\Gamma$ under the action of~$G$, where
two colorings $f$ and $f'$ are equivalent if $f' = f\circ g$ for some
$g\in G$.
\end{definition}

As a prerequisite for the following theorem, we adopt the notation $\Gamma/g$
from \cite{Cameron:2007}. For a graph $\Gamma$ and a permutation $g$ of its
vertex set, $\Gamma/g$ denotes the graph obtained by identifying the vertices
in each cycle of the disjoint-cycle decomposition of~$g$---that is, by
contracting each such cycle to a single vertex---and then replacing parallel
edges by single edges. This last step differs from the definition in
\cite{Cameron:2007}; however, since parallel edges do not affect the bivariate
chromatic polynomial, we suppress them for clarity.

\begin{theorem}
\label{thm:orbital}
For every graph $\Gamma$, every subgroup $G$ of $\Aut(\Gamma)$, every
$\lambda\in\mathbb{N}$ and $\mu=0,\ldots,\lambda$, we have
\begin{equation*}
OP_{\Gamma,G}(\lambda,\mu) \,=\, \frac{1}{|G|} \sum_{g\in G} P_{\Gamma/g}(\lambda,\mu).
\end{equation*}
\end{theorem}

\begin{proof}
We apply Burnside's lemma. Evidently, $(f,g)\mapsto f\circ g$ defines a right
group action of $G$ on the set of $\mu$-proper $\lambda$-colorings
of~$\Gamma$. For every $g\in G$, the fixpoints of $g$ under this action are
exactly the $\mu$-proper $\lambda$-colorings of $\Gamma$ that assign the same
color to all vertices within each cycle of~$g$. Since each of these colorings
corresponds uniquely to a $\mu$-proper $\lambda$-colorings of $\Gamma/g$, and
vice versa, the statement follows.
\end{proof}

The preceding theorem expresses $OP_{\Gamma,G}(\lambda,\mu)$ as an average of
bivariate chromatic polynomials in $\lambda$ and~$\mu$, and hence shows that
$OP_{\Gamma,G}(\lambda,\mu)$ is itself a polynomial in these variables. This
motivates the following definition.

\begin{definition}
For every graph $\Gamma$ and every subgroup $G$ of $\Aut(\Gamma)$ we refer to 
$OP_{\Gamma,G}(\lambda,\mu)$ as the \emph{orbital bivariate chromatic
  polynomial of $\Gamma$ with respect to $G$}.
\end{definition}

\begin{remark}
We highlight two fundamental specializations of the orbital bivariate
chromatic polynomial.
\begin{itemize}
\item If $\lambda=\mu$, then
\[
OP_{\Gamma,G}(\lambda,\lambda)=OP_{\Gamma,G}(\lambda),
\]
which is the orbital chromatic polynomial introduced in \cite{Cameron:2007}.
\item If $G=\{\textit{id}\}$, then
\[
OP_{\Gamma,\{\textit{id}\}}(\lambda,\mu)=P_\Gamma(\lambda,\mu),
\]
which is the bivariate chromatic polynomial introduced in~\cite{Dohmen:2003}.
\end{itemize}
Thus, the orbital bivariate chromatic polynomial simultaneously generalizes
both the orbital chromatic polynomial and the bivariate chromatic polynomial.
\end{remark}

\subsection{An illustrative example}

The following example continues Example~\ref{example:bivar-chrom-polyn-1} from
Subsection~\ref{sec:an-illustr-example}, now considering the orbital bivariate
chromatic polynomial of $\Gamma$ with respect to its automorphism group
$G = \Aut(\Gamma)$.

\begin{example}
Consider the graph $\Gamma = \exgra$, depicted in Figure~\ref{fig:example}.
Its automorphism group in cycle notation is
\begin{gather}
\label{eq:28}
\Aut(\Gamma) = \{\textit{id}, (13), (24), (13)(24)\}.
\end{gather}

\begin{table}[ht]
\centering
\begin{tabular}{c c c}
\hline
$g$ & $\Gamma/g$ & $P_{\Gamma/g}(\lambda,\mu)$ \\
\hline
$\textit{id}$ 
& $\Gamma/\textit{id} \cong \Gamma$ 
& $\lambda^4 - 5\lambda^2\mu +6\lambda\mu + 2\mu^2 -4\mu$ \\
$(13)$ 
& $\Gamma/(13) \cong \Kappa_3$ 
& $\lambda^3 - 3\lambda\mu +2\mu$ \\
$(24)$ 
& path of length~$2$ with a loop at its central vertex 
& $\lambda^2(\lambda-\mu)$ \\
$(13)(24)$ 
& path of length~$1$ with a loop at one end 
& $\lambda(\lambda-\mu)$ \\
\hline
\end{tabular}
\caption{Quotient graphs and their bivariate chromatic polynomials}
\label{tab:example-quotient-polynomials}
\end{table}

Table~\ref{tab:example-quotient-polynomials} shows the quotient graph
$\Gamma/g$ along with its bivariate chromatic polynomial for each
$g\in\Aut(\Gamma)$.  Taking the average of these polynomials, we obtain the
orbital bivariate chromatic polynomial of $\Gamma$ with respect to
$\Aut(\Gamma)$:
\begin{gather*}
  OP_{\Gamma,\Aut(\Gamma)}(\lambda,\mu) =
                                          \frac{1}{4} \left( \lambda^4 + 2 \lambda^3 - 6\lambda^2 \mu + \lambda^2 +
                                          2\lambda \mu + 2\mu^2 - 2\mu
                                          \right).
\end{gather*}                                          
For $\lambda=\mu$, this specializes to the orbital chromatic polynomial
\begin{gather*}
OP_{\Gamma,\Aut(\Gamma)}(\lambda) =
                                              \frac{1}{4} \lambda(\lambda-1)^2(\lambda-2),
                                              \end{gather*}
                                              which appears in
                                              \cite{Cameron:2008}.
\end{example}

\subsection{Fundamental properties}

This section establishes several fundamental properties of the orbital
bivariate chromatic polynomial. The first theorem concerns its total degree
and its degree with respect to~$\lambda$. We begin by proving a corresponding
result for the bivariate chromatic polynomial.

\begin{lemma}
\label{lem:bivariatedegree}
For every graph $\Gamma$,
$P_{\Gamma}(\lambda,\mu) = \lambda^{n(\Gamma)} + Q_\Gamma(\lambda,\mu)$, where
$Q_\Gamma(\lambda,\mu)\in \mathbb{Z}[\lambda,\mu]$ with
$\deg Q_\Gamma(\lambda,\mu)\le n(\Gamma)-1$.
In particular, $\deg P_{\Gamma}(\lambda,\mu) = n(\Gamma)$.
\end{lemma}

\begin{proof}
By \cite[Theorem~9]{Dohmen:2003}, the bivariate chromatic polynomial can be
written as
\begin{gather*}
P_{\Gamma}(\lambda,\mu) \,=\, \sum_{\substack{k,l=0\\ 0\le l\le k}}^{m} (-1)^k b_{k,l}
\lambda^{n(\Gamma)-k-l} \mu^l\, =\, \lambda^{n(\Gamma)} \,+  \! \sum_{\substack{k,l=0\\ 0\le l\le k>0}}^{m} (-1)^k b_{k,l}
\lambda^{n(\Gamma)-k-l} \mu^l ,
\end{gather*}
where $m\in\mathbb{N}_0$, $b_{0,0}=1$ and $b_{k,l}\in\mathbb{N}_0$ for
$0\le l\le k\le m$.
\end{proof}

The following theorem extends Lemma~\ref{lem:bivariatedegree} from the
bivariate chromatic polynomial to its orbital counterpart. When
$G=\{\textit{id}\}$, it reduces to Lemma~\ref{lem:bivariatedegree}.

\begin{theorem}
\label{thm:orbital2}
For every graph $\Gamma$ and every subgroup $G$ of its automorphism group,
\begin{gather*}
OP_{\Gamma,G}(\lambda,\mu) = \frac{1}{|G|} \left( \lambda^{n(\Gamma)} +
Q_{\Gamma,G}(\lambda,\mu) \right),
\end{gather*}
where $Q_{\Gamma,G}(\lambda,\mu)\in \mathbb{Z}[\lambda,\mu]$ satisfies
$\deg Q_{\Gamma,G}(\lambda,\mu) \le n(\Gamma)-1$. In particular, both the
degree of $OP_{\Gamma,G}(\lambda,\mu)$ and its partial degree with respect to
$\lambda$ are equal to $n(\Gamma)$, and its leading term is given by $\frac{1}{|G|}\lambda^{n(\Gamma)}$.
\end{theorem}

\begin{proof}
By Theorem~\ref{thm:orbital} and Lemma~\ref{lem:bivariatedegree}, we have
\begin{gather*}
OP_{\Gamma,G}(\lambda,\mu) 
= \frac{1}{|G|} \left( \lambda^{n(\Gamma)} +
Q_{\Gamma}(\lambda,\mu) + \sum_{\substack{g\in G\\ g\neq\textit{id}}}
P_{\Gamma/g}(\lambda,\mu) \right),
\end{gather*}
where $\deg Q_{\Gamma}(\lambda,\mu) \le n(\Gamma)-1$. For every
$g\neq \textit{id}$ at least two vertices of $\Gamma$ are identified in
$\Gamma/g$, which implies that
\begin{gather*}
\deg P_{\Gamma/g}(\lambda,\mu) = n(\Gamma/g) \le n(\Gamma)-1.
\end{gather*}
Thus, by setting
\[ Q_{\Gamma,G}(\lambda,\mu) = Q_{\Gamma}(\lambda,\mu) + \sum_{\substack{g\in
    G\\ g\neq\textit{id}}} P_{\Gamma/g}(\lambda,\mu), \]
it follows that $\deg Q_{\Gamma,G}(\lambda,\mu) \le n(\Gamma)-1$,
which proves the statement.
\end{proof}

As a consequence of Theorem~\ref{thm:orbital2}, the leading term of the
orbital bivariate chromatic polynomial determines the order of the graph and
the size of its automorphism group.

\begin{corollary}
\label{cor:distinctive}
Let $\Gamma$ and $\Gamma'$ be graphs such that
\[
OP_{\Gamma,\Aut(\Gamma)}(\lambda,\mu)
=
OP_{\Gamma',\Aut(\Gamma')}(\lambda,\mu).
\]
Then,
\[
n(\Gamma)=n(\Gamma')
\quad\text{and}\quad
|\Aut(\Gamma)|=|\Aut(\Gamma')|.
\]
\end{corollary}

\begin{proof}
By Theorem~\ref{thm:orbital2}, the leading term of
$OP_{\Gamma,\Aut(\Gamma)}(\lambda,\mu)$ equals
$1/|\Aut(\Gamma)|\lambda^{n(\Gamma)}$.  Since the polynomials of $\Gamma$ and
$\Gamma'$ coincide, they have the same leading term.
\end{proof}

As a consequence of Corollary~\ref{cor:distinctive}, the orbital bivariate
chromatic polynomial distinguishes asymmetric graphs from graphs with
nontrivial automorphisms: their orbital bivariate chromatic polynomials are
different. Notably, examples of non-isomorphic graphs with identical
bivariate chromatic polynomials are known; see \cite{Dohmen:2003}. This
naturally raises the question of whether non-isomorphic graphs with identical
orbital bivariate chromatic polynomials exist. A related open question is
whether graphs with the same orbital bivariate chromatic polynomial necessarily
have isomorphic automorphism groups.

As observed in \cite{Dohmen:2003}, the bivariate chromatic polynomial of a
disjoint union of graphs is the product of the bivariate chromatic polynomials
of its components. Under suitable assumptions, an analogous result holds for
the orbital bivariate chromatic polynomial.

\begin{theorem}
Let $\Gamma$ be the disjoint sum of non-isomorphic connected graphs
$\Gamma_1$ and $\Gamma_2$. Then,
\begin{gather}
\label{eq:21}
OP_{\Gamma, \Aut(\Gamma)}(\lambda,\mu) \,=\,
OP_{\Gamma_1,\Aut(\Gamma_1)}(\lambda,\mu)\, OP_{\Gamma_2,\Aut(\Gamma_2)}(\lambda,\mu). 
\end{gather}
\end{theorem}

\begin{proof}
Without loss of generality, we consider $\Aut(\Gamma_1)$ and $\Aut(\Gamma_2)$
as subgroups of $\Aut(\Gamma)$. Because $\Gamma_1$ and $\Gamma_2$ are
non-isomorphic and connected, $\Aut(\Gamma)$ can be viewed as the internal
direct product of $\Aut(\Gamma_1)$ and $\Aut(\Gamma_2)$.  Accordingly, each
$g\in\Aut(\Gamma)$ can be uniquely written (up to order) as $g=g_1g_2$ with
$g_1\in\Aut(\Gamma_1)$ and $g_2\in\Aut(\Gamma_2)$.  Hence, setting
$G_1 = \Aut(\Gamma_1)$ and $G_2 = \Aut(\Gamma_2)$, we have
$|G_1G_2| = |G_1| |G_2|$. By Theorem~\ref{thm:orbital},
\begin{align}
  \label{eq:26}
OP_{\Gamma,\Aut(\Gamma)}(\lambda,\mu) & = \frac{1}{|G_1||G_2|} \sum_{g_1g_2\in
                                        G_1 G_2}
                                        P_{\Gamma/g_1g_2}(\lambda,\mu) .
\end{align}  
Evidently, for every $g_1\in G_1$ and $g_2\in G_2$, the graph $\Gamma/g_1g_2$
is the disjoint sum of $\Gamma_1/g_1$ and $\Gamma_2/g_2$. Hence, by the
multiplicativity of the bivariate chromatic polynomial \cite{Dohmen:2003}, we
have
\[ P_{\Gamma/g_1g_2}(\lambda,\mu) = P_{\Gamma_1/g_1}(\lambda,\mu)
P_{\Gamma_2/g_2}(\lambda,\mu). \] Substituting this into Eq.~\eqref{eq:26}
yields
\begin{align*}
  OP_{\Gamma,\Aut(\Gamma)}(\lambda,\mu) & = \frac{1}{|G_1||G_2|}
                                          \sum_{g_1 g_2\in G_1 G_2}
                                          P_{\Gamma_1/g_1}(\lambda,\mu)
                                          P_{\Gamma_2/g_2}(\lambda,\mu) \\
                                        & = \frac{1}{|G_1|} \left( \sum_{g_1\in G_1} P_{\Gamma_1/g_1}(\lambda,\mu)
                                          \right)
                                          \frac{1}{|G_2|} \left( \sum_{g_2\in G_2} P_{\Gamma_2/g_2}(\lambda,\mu)
                                          \right),
\end{align*}
which equals $OP_{\Gamma_1,G_1}(\lambda,\mu) \,
OP_{\Gamma_2,G_2}(\lambda,\mu)$, as stated in Eq.~\eqref{eq:21}.
\end{proof}

We conclude this section with an alternative combinatorial interpretation of
the orbital bivariate chromatic polynomial in the case $\lambda=2$ and
$\mu=1$.

\begin{theorem}
For every graph $\Gamma$ and every subgroup $G$ of its automorphism group,
$OP_{\Gamma,G}(2,1)$ counts the equivalence classes of independent
vertex subsets of $\Gamma$ with respect to $G$, where $I$ and $J$ are
equivalent with respect to $G$ if $J=g(I)$ for some $g\in G$.
\end{theorem}

\begin{proof}
Since any $1$-proper $2$-coloring corresponds to an independent vertex subset
of $\Gamma$, and vice versa, the statement follows immediately from
Definition~\ref{def:orbital}.
\end{proof}

\begin{example}
We again consider the graph $\Gamma = \exgra$, depicted in
Figure~\ref{fig:example}.  Among the six independent vertex subsets,
$\emptyset$, $\{1\}$, $\{2\}$, $\{3\}$, $\{4\}$, and $\{1,3\}$, the sets
$\{1\}$ and $\{3\}$, as well as $\{2\}$ and $\{4\}$, are equivalent with
respect to $\Aut(\Gamma)$, given by Eq.~\eqref{eq:28}. Thus, there are
$OP_{\Gamma,\Aut(\Gamma)}(2,1) = 4$ independent vertex subsets such that no
two of them are equivalent.
\end{example}

\section{Special graph families}
\label{sec:spec-graph-class}

In this section, we derive expansions for the orbital bivariate chromatic
polynomial of edgeless graphs, complete graphs, complete bipartite graphs,
stars, paths, cycles, and wheels. Each graph family is treated in a separate
subsection. We use the corresponding formulas for the bivariate chromatic
polynomial from Subsection~\ref{sec:explicit-formulas} without restating them
explicitly.

Some technical statements in Subsections~\ref{sec:pin} and~\ref{sec:rg}
concerning the structure of $\Gamma/g$, attributed to Kim et
al. \cite{Kim:2014}, are not fully substantiated in that source. To ensure
mathematical rigor, we therefore provide complete proofs of these statements.

Our results for complete bipartite graphs, stars, paths, cycles, and wheels are
also new in the setting of the orbital chromatic polynomial, obtained by
specializing to $\lambda=\mu$.

\subsection{Edgeless graphs}

Recall from Subsection~\ref{sec:explicit-formulas} that $\Kappa_n$ denotes the complete
graph on $n$ vertices. The orbital bivariate chromatic polynomial of its
complement, $\overline{\Kappa_n}$, can be readily determined from its
definition, and coincides with the orbital chromatic polynomial of
$\overline{\Kappa_n}$ as given in \cite{Cameron:2008,Cameron:2007}.

\begin{theorem}
\label{thm:edgeless}
For every $n,\lambda\in\mathbb{N}$ and each $\mu=0,\dots,\lambda$,
\begin{gather*}
OP_{\overline{\Kappa_n},\Aut(\overline{\Kappa_n}})(\lambda,\mu) =
{n+\lambda-1 \choose n}.
\end{gather*}
\end{theorem}

\begin{proof}
Since every $\lambda$-coloring of $\overline{\Kappa_n}$ is $\mu$-proper, each
equivalence class of $\mu$-proper $\lambda$-colorings corresponds to a
combination with repetition of $n$ colors chosen from $\lambda$ available
colors. By elementary combinatorics, there are ${n+\lambda-1 \choose n}$ such
combinations, and hence the same number of equivalence classes.
\end{proof}

\subsection{Complete graphs}

\newcommand{\asn}[2]{\left[ #1 \atop #2 \right]_2}

Without loss of generality, we may assume that the vertex set and edge set of
$\Kappa_n$ is given by $V(\Kappa_n)=\{1,\dots,n\}$ and
$E(\Kappa_n) = \{\{v,w\}\mid 1\le v < w \le n\}$, respectively.

The following lemma clarifies the structure of $\Kappa_n/g$ for
$g\in\Aut(\Kappa_n)$.  Clearly, $\Aut(\Kappa_n)=S_n$, where $S_n$ denotes the
symmetry group of $\{1,\dots,n\}$. For each $\sigma\in S_n$, we denote by
$|\sigma|$ the number of cycles in the disjoint cycle decomposition
of~$\sigma$, and by $|\sigma|_1$ the number of cycles of length one, i.e., the
number of fixpoints of $\sigma$.

\begin{lemma}
\label{lem:complete-graphs}
For every $n\in\mathbb{N}$ and every $\sigma\in S_n$, the graph
$\Kappa_n/\sigma$ is isomorphic to the graph obtained from $\Kappa_{|\sigma|}$ by
attaching loops to $|\sigma|-|\sigma|_1$ of its vertices.
\end{lemma}

\begin{proof}
Let $\sigma_1\dots \sigma_{|\sigma|}$ denote the disjoint cycle decomposition
of~$\sigma$. Identifying the vertices of $\Kappa_n$ within each cycle
$\sigma_i$ ($i=1,\dots,|\sigma|$) produces a graph on $|\sigma|$ vertices, in
which distinct vertices are adjacent, and vertices corresponding to cycles of
length greater than 1 carry a loop.
\end{proof}

The orbital bivariate chromatic polynomial of $\Kappa_n$ with respect to
$\Aut(\Kappa_n)$ is given by the following theorem.  Here, $\asn{n}{k}$
denotes the 2-associated Stirling number of the first kind, which counts the
number of permutations of $n$ elements that decompose into exactly $k$ cycles,
each cycle having length at least 2 (see \cite[p.~256]{Comtet:1974} and entry
A008306 in the OEIS~\cite{OEIS:2026}).

\begin{theorem}
\label{thm:complete-graphs-1}
For every $n,\lambda\in\mathbb{N}$ and each $\mu=0,\dots,\lambda$,
\begin{gather*}
OP_{\Kappa_n,\Aut(\Kappa_n)}(\lambda,\mu) = \frac{1}{n!} \sum_{m=0}^n
{n\choose m}  Q_m(\lambda-\mu) P_{\Kappa_{n-m}}(\lambda,\mu),
\end{gather*}
where $Q_m(x)\in\mathbb{Z}[x]$ is defined by
\begin{gather}
\label{eq:22}
Q_m(x) = \sum_{k=0}^m \asn{m}{k} x^k ,
\end{gather}
and $P_{\Kappa_{n-m}}(\lambda,\mu)$ is given by Eq.~\eqref{eq:8}.
\end{theorem}

\begin{proof}
By Lemma~\ref{lem:complete-graphs}, for any $\sigma\in S_n$ we have
\[ P_{\Kappa_n/\sigma}(\lambda,\mu) =
P_{\Kappa_{|\sigma|_1}}(\lambda,\mu)(\lambda-\mu)^{|\sigma|-|\sigma_1|} . \]
Therefore, by Theorem \ref{thm:orbital},
\begin{gather*}
OP_{\Kappa_n,\Aut(\Kappa_n)}(\lambda,\mu) = \frac{1}{|S_n|} \sum_{\sigma\in S_n}
  P_{\Kappa_n/\sigma}(\lambda,\mu) 
 = \frac{1}{n!} \sum_{m=0}^n P_{\Kappa_m}(\lambda,\mu) \sum_{\sigma\in S_n
                                            \atop |\sigma|_1=m} (\lambda-\mu)^{|\sigma|-m} .
\end{gather*}
Since there are ${n \choose m}$ ways to select $m$ fixpoints from
$\{1,\dots,n\}$, we find
\begin{gather*}
OP_{\Kappa_n,\Aut(\Kappa_n)}(\lambda,\mu)  
  = \frac{1}{n!} \sum_{m=0}^n  {n\choose m} P_{\Kappa_m}(\lambda,\mu) 
    \sum_{\sigma\in S_{n-m} \atop |\sigma|_1=0}
    (\lambda-\mu)^{|\sigma|}.
\end{gather*}
By symmetry of the binomial coefficients, this can be written as
\begin{gather*}
OP_{\Kappa_n,\Aut(\Kappa_n)}(\lambda,\mu)  = \frac{1}{n!} \sum_{m=0}^n  {n\choose m} P_{\Kappa_{n-m}}(\lambda,\mu) 
    \sum_{\sigma\in S_{m} \atop |\sigma|_1=0}
    (\lambda-\mu)^{|\sigma|} . 
\end{gather*}
Finally, by Eq.~\eqref{eq:22}, the inner sum simplifies to $Q_m(\lambda-\mu)$,
which completes the proof.
\end{proof}

\begin{remark}
By the inclusion-exclusion principle, for any $m\in\mathbb{N}_0$ and
$k=0,\dots,m$, we have
\begin{gather*}
\asn{m}{k} = \sum_{j=0}^k (-1)^j {m\choose j}\left[ m-j \atop
k-j\right],
\end{gather*}
where the bracketed term on the right-hand side denotes an unsigned Stirling
number of the first kind. 
\end{remark}

For $n=1,\dots,5$, the orbital bivariate chromatic polynomials given by
Theorem~\ref{thm:complete-graphs-1} are shown in
Table~\vref{tab:orbital_complete}.

\renewcommand{\arraystretch}{1.2}
\begin{table}[!t]
\footnotesize
\begin{center}
\begin{tabular}{rc}
  $n$ & $OP_{\Kappa_n,\Aut(\Kappa_n)}(\lambda,\mu)$ \\[1pt] \hline
1 & $\lambda$ \\
2 & $\frac{1}{2} \, {\lambda}^{2} + \frac{1}{2} \, {\lambda} - {\mu}$ \\
3 & $\frac{1}{6} \, {\lambda}^{3} + \frac{1}{2} \, {\lambda}^{2} - {\lambda}
    {\mu} + \frac{1}{3} \, {\lambda}$ \\
4 & $\frac{1}{24} \, {\lambda}^{4} + \frac{1}{4} \, {\lambda}^{3} - \frac{1}{2}
    \, {\lambda}^{2} {\mu} + \frac{11}{24} \, {\lambda}^{2} - \frac{1}{2} \,
    {\lambda} {\mu} + \frac{1}{2} \, {\mu}^{2} + \frac{1}{4} \, {\lambda} -
    \frac{1}{2} \, {\mu}$ \\
5 & $\frac{1}{120} \, {\lambda}^{5} + \frac{1}{12} \, {\lambda}^{4} -
    \frac{1}{6} \, {\lambda}^{3} {\mu} + \frac{7}{24} \, {\lambda}^{3} -
    \frac{1}{2} \, {\lambda}^{2} {\mu} + \frac{1}{2} \, {\lambda} {\mu}^{2} +
    \frac{5}{12} \, {\lambda}^{2} - \frac{5}{6} \, {\lambda} {\mu} +
    \frac{1}{5} \, {\lambda}$
  \\[1pt] \hline
\end{tabular}
\caption{Orbital bivariate chromatic polynomials of $\Kappa_n$ with respect to
  $\Aut(\Kappa_n)$}
\label{tab:orbital_complete}
\end{center}
\end{table}

\subsection{Complete bipartite graphs}

For the complete bipartite graph $\Kappa_{m,n}$, we assume
$V(\Kappa_{m,n}) = \{1,\dots,m+n\}$ and
$E(\Kappa_{m,n}) = \{\{i,j\}\mid 1\le i\le m < j \le m+n\}$, and define
\begin{align}
 A_{m,n} & = \left\{\sigma\tau \mid \sigma,\tau\in S_{m+n},
\{1,\dots,m\}\subseteq \Fix(\sigma), \{m+1,\dots,m+n\}\subseteq \Fix(\tau)
           \right\}, \label{sigmatau}\\
B_{n,n} & = \{\sigma\in S_{2n}\mid
          \text{$\sigma(i)>n$, for $i=1,\dots,n$}\}, \label{bnn}
\end{align}
where $\Fix(\sigma)$ denotes the set of fixpoints of~$\sigma$ (similarly for
$\tau$).  Clearly,
\begin{align}
\label{eq:15}   
\Aut(\Kappa_{m,n}) & =
\begin{cases}
A_{m,n}, &  \text{if $m\neq n$}, \\ 
  A_{n,n} \cup B_{n,n}, & \text{if $m=n$};
\end{cases} \\
\intertext{and}
\label{eq:16}        
|\Aut(\Kappa_{m,n})| & =
\begin{cases}
m!n!, &  \text{if $m\neq n$}, \\ 
2n!n!, & \text{if $m=n$}.
\end{cases}               
\end{align}

\begin{lemma}
\label{lem:compl-bipart-graphs}
For every $m,n\in\mathbb{N}$ and every $\sigma\tau\in A_{m,n}$, with
$\sigma,\tau$ as in Eq.~\eqref{sigmatau}, the graph $\Kappa_{m,n}/\sigma\tau$
is isomorphic to the complete bipartite graph $\Kappa_{|\sigma'|,|\tau'|}$,
where $|\sigma'|$ and $|\tau'|$ denote the number of cycles in the disjoint
cycle decomposition of the restricted permutations
$\sigma'=\sigma|_{\{m+1,\dots,m+n\}}$ and $\tau'=\tau|_{\{1,\dots,m\}}$,
respectively.  Moreover, for every $\sigma\in B_{n,n}$, the graph
$\Kappa_{n,n}/\sigma$ is a complete graph on $|\sigma|$ vertices, with loops
attached to all vertices.
\end{lemma}

\begin{proof}
Let $\sigma\tau\in A_{m,n}$, with $\sigma,\tau$ as in Eq.~\eqref{sigmatau}.
Since $\sigma$ fixes all of $\{1,\dots,m\}$ and $\tau$ fixes all of
$\{m+1,\dots,m+n\}$, we have $\sigma\tau = \sigma'\tau'$.  Hence, the disjoint
cycle decomposition of $\sigma\tau$ can be written as
$\sigma_1'\dots\sigma_{|\sigma'|}' \tau_1'\dots\tau_{|\tau'|}'$ where
$\sigma_1'\dots\sigma_{|\sigma'|}'$ is the disjoint cycle decomposition of
$\sigma'$, and $\tau_1'\dots\tau_{|\tau'|}'$ is the disjoint cycle
decomposition of $\tau'$.  Identifying the vertices of $\Kappa_n$ within each
cycle $\sigma_i'$ ($i=1,\dots,|\sigma'|$) and within each cycle $\tau_j'$
($j=1,\dots,|\tau'|$) produces a graph on $|\sigma'| + |\tau'|$ vertices, in
which precisely the vertices corresponding to $\sigma_i'$ and $\tau_j'$ are
adjacent, for $i=1,\dots,|\sigma'|$ and $j=1,\dots,|\tau'|$. This proves the
first statement of the lemma.
\par
For the second statement, consider $\sigma\in B_{n,n}$, with $\sigma$ as in
Eq.~\eqref{bnn}. Let $\sigma_1\dots\sigma_{|\sigma|}$ denote the disjoint
cycle decomposition of~$\sigma$. Each cycle $\sigma_i$ consists of an
alternating sequence of vertices from $\{1,\dots,n\}$ and of vertices from
$\{n+1,\dots,2n\}$. Identifying the vertices of $\Kappa_{n,n}$ within each
cycle $\sigma_i$ ($i=1,\dots,|\sigma|$), yields a complete graph on $|\sigma|$
vertices, in which each vertex is incident with a loop.
\end{proof}

Recall that we use the bracketed notation $\left[ n \atop k \right]$ for the
unsigned Stirling numbers of the first kind.

\begin{theorem}
\label{thm:compl-bipart-graphs}
For every $m,n,\lambda\in\mathbb{N}$ and each $\mu=0,\dots,\lambda$,
\begin{align}
  \label{eq:17}
  OP_{\Kappa_{m,n},\Aut(\Kappa_{m,n})}(\lambda,\mu) & =
  \frac{1}{m!n!}\, \sum_{k=0}^m\, \sum_{l=0}^n \, \left[ m\atop k\right] \left[n\atop
                                                      l\right]
                                                      P_{\Kappa_{k,l}}(\lambda,\mu)
                                                      \quad (m\neq n);
  \\
  \label{eq:18}
  OP_{\Kappa_{n,n},\Aut(\Kappa_{n,n})}(\lambda,\mu) & =
  \frac{1}{2n!n!}  \sum_{k,l=0}^n \left[ n\atop k\right] \left[n\atop
                                                      l\right]
                                                      P_{\Kappa_{k,l}}(\lambda,\mu)
                                                      \,+\,\frac{1}{2}
                                                      {n+\lambda-\mu-1 \choose n},
\end{align}
where $P_{\Kappa_{k,l}}(\lambda,\mu)$ is given by Eq.~\eqref{eq:19}.
\end{theorem}

\begin{proof}
By Lemma~\ref{lem:compl-bipart-graphs} and the definition of $A_{m,n}$
in Eq.~\eqref{sigmatau}, we have
\begin{gather*}
\sum_{\sigma\tau\in A_{m,n}} P_{\Kappa_{m,n}/\sigma\tau}(\lambda,\mu)  =
\sum_{\sigma\tau\in A_{m,n}}
P_{\Kappa_{|\sigma'|,|\tau'|}}(\lambda,\mu) 
= \sum_{\sigma'\in S_n} \sum_{\tau'\in S_m}
P_{\Kappa_{|\sigma'|,|\tau'|}}(\lambda,\mu),
\end{gather*}
where the first and second sum is over all $\sigma\tau\in A_{m,n}$, with
$\sigma,\tau$ as in Eq.~\eqref{sigmatau}. Thus, we obtain
\begin{gather} 
\label{eq:10}
\sum_{\sigma\tau\in A_{m,n}} P_{\Kappa_{m,n}/\sigma\tau}(\lambda,\mu)  =
\sum_{k=0}^m \sum_{l=0}^n \left[ m \atop k \right] \left[ n \atop l \right]
P_{\Kappa_{k,l}}(\lambda,\mu) .
\end{gather}
From Lemma~\ref{lem:compl-bipart-graphs} and the definition of $B_{n,n}$
in Eq.~\eqref{bnn}, we obtain
\begin{align}  
  \sum_{\sigma\in B_{n,n}} P_{\Kappa_{n,n}/\sigma}(\lambda,\mu)   & =
                                                                    \sum_{\sigma\in B_{n,n}} (\lambda-\mu)^{|\sigma|} = 
                                                                    \sum_{k=0}^n
                                                                    \left[ n \atop k
                                                                    \right] n! \,
                                                                    (\lambda-\mu)^k
                                                                    = n!n! {n+\lambda-\mu-1 \choose n}, \label{eq:14}
\end{align}
where, in the last step, we used the identity (cf.~\cite{Graham:1994}),
\begin{gather}
\label{eq:20}
\frac{1}{n!} \sum_{k=0}^n \left[{n \atop k}\right] x^k = {n+x-1\choose n}.
\end{gather}
Finally, Eqs.~\eqref{eq:17} and \eqref{eq:18} follow from
Theorem~\ref{thm:orbital} and Eqs.~\eqref{eq:15}--\eqref{eq:16}
and \eqref{eq:10}--\eqref{eq:14}.
\end{proof}

\subsection{Stars}

The following theorem states the implication of
Theorem~\ref{thm:compl-bipart-graphs} for the star $\Kappa_{1,n}$.  We use
$x^{\overline{k}}$ to denote the $k$-th rising factorial of~$x$; that is,
$x^{\overline{k}} = x(x+1)\cdots(x+k-1)$.

\begin{theorem}
\label{thm:orbit-bivar-chrom}
For every $\lambda\in\mathbb{N}$ and each $\mu=0,\dots,\lambda$,
\begin{gather}
\label{eq:12}
OP_{\stern{1},\Aut(\stern{1})}(\lambda,\mu) = \frac{1}{2}\lambda^2 +
\frac{1}{2}\lambda - \mu\, .
\end{gather}
Moreover, for $n\ge 2$,
\begin{gather}
\label{eq:11}
OP_{\stern{n},\Aut(\stern{n})}(\lambda,\mu) =  \frac{1}{n!}\big(\lambda^2 +
(n-1)\lambda - n\mu\big)\lambda^{\overline{n-1}} \, .
\end{gather}
\end{theorem}

\begin{proof}
Eq.~\eqref{eq:12} follows from Theorem~\ref{thm:compl-bipart-graphs} with
$m=n=1$.  To prove Eq.~\eqref{eq:11}, we apply
Theorem~\ref{thm:compl-bipart-graphs} to $m=1$ and $n\ge 2$, and then
use Eqs.~\eqref{eq:3} and \eqref{eq:20}.  Thus, we obtain
\begin{align}
  OP_{\stern{n},\Aut(\stern{n})} (\lambda,\mu) & = \frac{1}{n!} \sum_{l=0}^n \left[ n \atop l\right]
                                               P_{\stern{l}}(\lambda,\mu)\, . \notag \\                               
 & = \frac{1}{n!} \sum_{l=0}^n \left[ n
                                               \atop l \right] \left( \lambda^l(\lambda-\mu) + (\lambda-1)^l\mu
                                               \right) \notag \notag \\
                                             & =  {n+\lambda-1 \choose n} (\lambda-\mu) + {n+\lambda-2\choose
                                               n}
                                               \mu\,   \label{binom} \\
& = \frac{\lambda^{\overline{n}}}{n!} \, (\lambda-\mu) +
  \frac{(\lambda-1)\lambda^{\overline{n-1}}}{n!} \, \mu \notag \\
& = \frac{1}{n!}\big( (\lambda+n-1)(\lambda-\mu) + (\lambda-1)\mu \big)
  \lambda^{\overline{n-1}}. \notag
\end{align}
Grouping the terms in parentheses according to powers of $\lambda$ yields the
desired result.
\end{proof}

\begin{remark}
Eq.~\eqref{binom}, and hence Theorem~\ref{thm:orbit-bivar-chrom}, can also be
deduced from Theorem~\ref{thm:edgeless}, since
\[ OP_{\Kappa_{1,n},\Aut(\Kappa_{1,n})} (\lambda,\mu) =
OP_{\overline{\Kappa_n},\Aut(\overline{\Kappa_n})}(\lambda,\mu)\,(\lambda-\mu)
+ OP_{\overline{\Kappa_n},\Aut(\overline{\Kappa_n})}(\lambda-1,\mu)\,\mu . \]
This equation follows from a case distinction for the central vertex of the
star: if it is colored with one of the $\lambda - \mu$ colors greater than
$\mu$, then its $n$ neighbors are subject to no restrictions. If, however, it
is colored with one of the $\mu$ colors less than or equal to $\mu$, then only
$\lambda - 1$ colors remain available for its neighborhood.
\end{remark}

For $n=1,\dots,6$, the orbital bivariate chromatic polynomials given by
Theorem~\ref{thm:orbit-bivar-chrom} are shown in factored form in
Table~\vref{tab:orbital_star_Aut}.

\renewcommand{\arraystretch}{1.2}
\begin{table}[!t]
\footnotesize
\begin{center}
\begin{tabular}{rc}
  $n$ & $OP_{\stern{n},\Aut(\stern{n})}(\lambda,\mu)$ \\[1pt] \hline
1 & $\frac{1}{2} \, \lambda^{2} + \frac{1}{2} \, \lambda - \mu$ \\
2 & $\frac{1}{2} \, {\left({\lambda}^{2} + {\lambda} - 2 \, {\mu}\right)}
    {\lambda}$ \\
3 & $\frac{1}{6} \, {\left({\lambda}^{2} + 2 \, {\lambda} - 3 \, {\mu}\right)}
    {\left({\lambda} + 1\right)} {\lambda}$ \\
4 & $\frac{1}{24} \, {\left({\lambda}^{2} + 3 \, {\lambda} - 4 \,
    {\mu}\right)} {\left({\lambda} + 2\right)} {\left({\lambda} + 1\right)}
    {\lambda}$ \\
5 & $\frac{1}{120} \, {\left({\lambda}^{2} + 4 \, {\lambda} - 5 \,
    {\mu}\right)} {\left({\lambda} + 3\right)} {\left({\lambda} + 2\right)}
    {\left({\lambda} + 1\right)} {\lambda}$ \\
6 & $\frac{1}{720} \, {\left({\lambda}^{2} + 5 \, {\lambda} - 6 \, {\mu}\right)} {\left({\lambda} + 4\right)} {\left({\lambda} + 3\right)} {\left({\lambda} + 2\right)} {\left({\lambda} + 1\right)} {\lambda}$
\\[1pt] \hline
\end{tabular}
\caption{Orbital bivariate chromatic polynomials of $\stern{n}$ with respect
  to $\Aut(\stern{n})$}
\label{tab:orbital_star_Aut}
\end{center}
\end{table}

\subsection{Paths}
\label{sec:pin}

For the path $\Pi_n$ of length $n$, we henceforth assume that its vertex set
and edge set are $V(\Pi_n)=\{0,\dots,n\}$ and
$E(\Pi_n) = \{\{v,v+1\} \mid v=0,\dots,n-1\}$.

The automorphism group of $\Pi_n$ is easy to describe: for $n\ge 1$,
$\Aut(\Pi_n)=\{\textit{id},\pi\}$, where \textit{id} is the identity and
$\pi(v)=n-v$ for $v=0,\dots,n$. For $n=0$, we have
$\Aut(\Pi_0)=\{\textit{id}\}$.

\begin{lemma}[\cite{Kim:2014}]
\label{lem:reflectpath}
For every $n\in\mathbb{N}_0$, the graph $\Pi_n/\pi$ is a path of length $\lfloor
n/2\rfloor$, with a loop attached to one of its end vertices if $n$ is odd.
\end{lemma}

\begin{proof}
Evidently, the disjoint cycle decomposition of $\pi$ is 
\[ \pi = (0,n)(1,n-1)\dots (n/2-1,n/2+1)(n/2) \]
if $n$ is even, and
\[ \pi = (0,n)(1,n-1)\dots ((n-1)/2,(n+1)/2) \] if $n$ is odd.  Following the
construction of $\Pi_n/\pi$, we identify the vertices within each cycle.  This
gives a path of length $n/2$ if $n$ is even, and a path of length $(n-1)/2$
with a loop attached to the vertex corresponding to the cycle
$((n-1)/2,(n+1)/2)$ if $n$ is odd.
\end{proof}

The following theorem provides closed-form expansions for the orbital
bivariate chromatic polynomial of $\Pi_n$ with respect to its automorphism
group.  The respective bivariate chromatic polynomials are given in
Eq.~\eqref{eq:6}.

\begin{theorem}
\label{thm:path}
For every $n\in\mathbb{N}_0$, every $\lambda\in\mathbb{N}$, and each
$\mu=0,\dots,\lambda$, we have
\begin{align}
\label{eq:2}
OP_{\Pi_n,\Aut(\Pi_n)}(\lambda,\mu) & = \frac{1}{2}\left(
                                      P_{\Pi_n}(\lambda,\mu) +
                                      P_{\Pi_{n/2}}(\lambda,\mu) \right), \\
  \intertext{if $n$ is even, respectively}
  \label{eq:9}
OP_{\Pi_n,\Aut(\Pi_n)}(\lambda,\mu) & = \frac{1}{2}\left(
                                      P_{\Pi_n}(\lambda,\mu) + (\lambda-\mu)
                                      P_{\Pi_{(n-3)/2}}(\lambda,\mu) \right),                                       
\end{align}
if $n$ is odd, where $\Pi_{-1}$ is considered as the empty graph.
\end{theorem}

\begin{proof}
By Theorem~\ref{thm:orbital} and $\Aut(\Pi_n)=\{\textit{id},\pi\}$,
\begin{equation}
\label{eq:7}
OP_{\Pi_n,\Aut(\Pi_n)}(\lambda,\mu) = \frac{1}{2} \left(
P_{\Pi_n/\textit{id}}(\lambda,\mu) + P_{\Pi_n/\pi}(\lambda,\mu)\right) ,
\end{equation}
where, trivially, $\Pi_n/\textit{id} = \Pi_n$.
For the second term in Eq.~(\ref{eq:7}), we distinguish whether $n$ is even or odd. If $n$ is
even, then by Lemma~\ref{lem:reflectpath}, $\Pi_n/\pi = \Pi_{n/2}$,
which implies Eq.~(\ref{eq:2}). If $n$ is odd, then by
Lemma~\ref{lem:reflectpath},
$\Pi_n/\pi$ is a path of length $(n-1)/2$ with a loop attached to one of its
end vertices. For the color of this end vertex there are $\lambda-\mu$
choices, while for the remaining vertices there are
$P_{\Pi_{(n-3)/2}}(\lambda,\mu)$ choices. Thus, for the second term in
Eq.~(\ref{eq:7}), we conclude that 
$P_{\Pi_n/\pi}(\lambda,\mu) = (\lambda-\mu) P_{\Pi_{(n-3)/2}}(\lambda,\mu)$,
which finally proves Eq.~(\ref{eq:9}).
\end{proof}

For $n=0,\dots,6$, the orbital bivariate chromatic polynomials given by
Theorem~\ref{thm:path} are shown in Table~\vref{tab:orbital_path_Aut}.

\renewcommand{\arraystretch}{1.2}
\begin{table}[!t]
\footnotesize
\begin{center}
\begin{tabular}{rc}
  $n$ & $OP_{\Pi_n,\Aut(\Pi_n)}(\lambda,\mu)$ \\[1pt] \hline
0 & $\lambda$ \\
1 & $\frac{1}{2} \, \lambda^{2} + \frac{1}{2} \, \lambda - \mu$ \\
2 & $\frac{1}{2} \, {\left(\lambda^{2} + \lambda - 2 \, \mu\right)} \lambda$ \\
3 & $\frac{1}{2} \, \lambda^{4} - \frac{3}{2} \, \lambda^{2} \mu + \frac{1}{2} \, \lambda^{2} +
    \frac{1}{2} \, \lambda \mu + \frac{1}{2} \, \mu^{2} - \frac{1}{2} \, \mu$ \\
4 & $\frac{1}{2} \, \lambda^{5} - 2 \, \lambda^{3} \mu + \frac{1}{2} \, \lambda^{3} + \frac{3}{2}
    \, \lambda^{2} \mu + \frac{3}{2} \, \lambda \mu^{2} - 2 \, \lambda \mu - \mu^{2} + \mu$ \\
5 & $\frac{1}{2} \, \lambda^{6} - \frac{5}{2} \, \lambda^{4} \mu + 2 \, \lambda^{3} \mu + 3 \, \lambda^{2}
    \mu^{2} + \frac{1}{2} \, \lambda^{3} - 2 \, \lambda^{2} \mu - 3 \, \lambda \mu^{2} - \frac{1}{2}
    \, \mu^{3} + \frac{1}{2} \, \lambda \mu + 2 \, \mu^{2} - \frac{1}{2} \, \mu$ \\
6 & $\frac{1}{2} \, \lambda^{7} - 3 \, \lambda^{5} \mu + \frac{5}{2} \, \lambda^{4} \mu + 5 \, \lambda^{3}
    \mu^{2} + \frac{1}{2} \, \lambda^{4} - 2 \, \lambda^{3} \mu - 6 \, \lambda^{2} \mu^{2} - 2 \, \lambda
    \mu^{3} + \frac{9}{2} \, \lambda \mu^{2} + \frac{3}{2} \, \mu^{3} - \frac{3}{2} \,
    \mu^{2}$ \\[1pt] \hline
\end{tabular}
\caption{Orbital bivariate chromatic polynomials of $\Pi_n$ with respect to
  $\Aut(\Pi_n)$}
\label{tab:orbital_path_Aut}
\end{center}
\end{table}

\subsection{Cycles}
\label{sec:rg}

For the cycle $\Gamma_n$, we set $V(\Gamma_n)=\{0,\dots,n-1\}$ and
$E(\Gamma_n)=\{\{v,(v+1) \bmod n\} \mid v=0,\dots,n-1\}$, where $n\ge 1$.  In
the special case $n=1$ resp.\ $n=2$, the cycle $\Gamma_n$ consists of a single
vertex with a loop attached, respectively of two vertices joined by parallel
edges.

The automorphism group of $\Gamma_n$ consists of $n$ rotations and $n$
reflections, which for $n\ge 3$ is known as the \emph{Dihedral group} of
order~$2n$.  Depending on whether $n$ is odd or even,
\begin{gather}
\AG{n}=\{r_0,\dots,r_{n-1},s_0,\dots,s_{n-1}\} \quad \text{($n$
  odd)}, \label{ag1} \\
\intertext{respectively}
\AG{n}=\{r_0,\dots,r_{n-1},s_0,\dots,s_{n/2-1,}, s_0',\dots,s_{n/2-1}'\}
\quad \text{($n$ even)}, \label{ag2}
\end{gather}
where, in both cases, $r_m(v)$, $s_m(v)$ and
$s_m'(v)$ for $v=0,\dots,n-1$ are given by
\begin{align*}
  r_m(v) & =  (v+m) \bmod n\, , \\
  s_m(v) & =  (2m-v) \bmod n\, , \notag  \\
  s_m'(v) & = (2m+1-v) \bmod n \, . \notag
\end{align*}
An important subgroup of $\AG{n}$ is its subgroup of rotations,
\begin{equation}
\label{eq:4}
\RG{n} = \{r_0,\dots,r_{n-1}\},
\end{equation}
which for $n=1$ and $n=2$ coincides with $\AG{n}$.

\begin{lemma}[\cite{Kim:2014}]
\label{lem:rotate}
For every $n\in\mathbb{N}$ and each $m=0,\dots,n-1$, the graph
$\Gamma_n/r_m$ is
\begin{enumerate}[(a)]
\item a cycle of length $\gcd(m,n)$ if $\gcd(m,n)\neq 2$;
\item a path of length 1 if $\gcd(m,n)=2$.
\end{enumerate}
\end{lemma}

\begin{proof}
Let $\sigma_0 \sigma_1 \dots \sigma_{k-1}$ be the disjoint cycle decomposition
of $r_m$. Without loss of generality, we may assume $i\in \sigma_i$ for
$i=0,\dots,k-1$.  Following the construction of $\Gamma_n/r_m$ we identify
vertices $x$ and $y$ of $\Gamma_n$ if $x$ and $y$ belong to the same cycle
$\sigma_i$; that is, if $x=r_m^s(i)=(i+sm)\bmod{n}$ and
$y=r_m^t(i)=(i+tm)\pmod{n}$ for some $i$, $s$ and $t$, or equivalently, if
there is a simultaneous solution to $i \equiv x \pmod{m}$ and
$i \equiv y \pmod{n}$.  By the generalized Chinese Remainder Theorem, this is
the case if and only if $x\equiv y \pmod{\gcd(m,n)}$. Thus, each cycle
$\sigma_i$ consists of vertices which are congruent to $i$ modulo
$k=\gcd(m,n)$. Each such cycle defines a vertex $\overline{\sigma_i}$ in
$\Gamma_n/r_m$, and any two (not necessarily distinct) vertices
$\overline{\sigma_i}$ and $\overline{\sigma_j}$ are joined by an edge in
$\Gamma_n/r_m$ if there exist $v\in \sigma_i$ and $w\in\sigma_j$ such that $v$
and $w$ are adjacent in $\Gamma_n$; that is, $v\equiv i \pmod{k}$ and
$w\equiv j \pmod{k}$ for some $v,w\in\{0,\dots,n-1\}$ such that
$v\equiv w\pm 1\pmod{n}$, which implies $i\equiv j\pm 1 \pmod{k}$ since
$k\mid n$.  On the other hand, if $i\equiv j\pm 1 \pmod{k}$ we show that
$\overline{\sigma_i}$ and $\overline{\sigma_j}$ are joined by an edge in
$\Gamma_n/r_m$. Without loss of generality we may assume that
$i\equiv j+1 \pmod{k}$, otherwise we exchange $i$ and $j$. We distinguish two
cases:

\emph{Case 1:} If $i>0$, then $i=j+1$. Since $i$ and $j$ are adjacent in
$\Gamma_n$, $\overline{\sigma_i}$ and $\overline{\sigma_j}$ are joined by an
edge in $\Gamma_n/r_m$.

\emph{Case 2:} If $i=0$, then $j=k-1\equiv n-1 \pmod{k}$ since $k\mid n$.
Therefore, $0\in\sigma_i$ and $n-1\in\sigma_j$. Since $0$ and $n-1$ are
adjacent in $\Gamma_n$, $\overline{\sigma_i}$ and $\overline{\sigma_j}$ are
adjacent in $\Gamma_n/r_m$.

Hence, for $k=\gcd(m,n)$, $\Gamma_n/r_m$ consists of the cycle
$(\overline{\sigma_0},\overline{\sigma_1},\dots,\overline{\sigma_{k-1}},\overline{\sigma_0})$
if $k\neq 2$ (which is a loop on $\overline{\sigma_0}$ if $k=1$) and of the path $(\overline{\sigma_0},\overline{\sigma_1})$ if $k=2$.
\end{proof}

\begin{lemma}[\cite{Kim:2014}]
\label{lem:reflect1}
For every $n\in\mathbb{N}$ and each $m=0,\ldots,n-1$, the graph
$\Gamma_n/s_m$ is a path of length $\lfloor n/2 \rfloor$, with a loop attached
to one of its end vertices if $n$ is odd.
\end{lemma}

\begin{proof}
Let $\sigma_0 \sigma_1 \ldots \sigma_{\lfloor n/2 \rfloor}$ be the disjoint
cycle decomposition of $s_m$, where
\begin{equation*}
\sigma_i = ((m-i) \bmod n, (m+i) \bmod n) \quad (i=0,\dots,\lfloor n/2 \rfloor).
\end{equation*}
Each cycle $\sigma_i$ defines a vertex $\overline{\sigma_i}$ in
$\Gamma_n/s_m$, and any two (not necessarily distinct) vertices
$\overline{\sigma_i}$ and $\overline{\sigma_j}$ are adjacent in $\Gamma_n/s_m$
if there exist adjacent vertices $v,w$ in $\Gamma_n$ such that $v\in\sigma_i$
and $w\in\sigma_j$; that is, $v\equiv w\pm 1 \pmod{n}$,
$v = (m\pm i) \bmod{n}$, and $w = (m\pm j) \bmod{n}$.  The conjunction of
these three conditions is equivalent to $i = j\pm 1$ or $i=j=(n-1)/2$, where
for the second alternative $n$ is required to be odd. Therefore,
$\Gamma_n/s_m$ is a path
$(\overline{\sigma_0},\overline{\sigma_1},\dots,\overline{\sigma_{\lfloor n/2
    \rfloor}})$ with an additional loop at
$\overline{\sigma_{\lfloor n/2 \rfloor}}$ in case that $n$ is odd.
\end{proof}

\begin{lemma}[\cite{Kim:2014}]
\label{lem:reflect2}
For every even $n\in\mathbb{N}$ and each $m=0,\dots,\frac{n}{2}-1$, the graph
$\Gamma_n/s_m'$ is a path of length $\frac{n}{2}-1$ with a loop attached to
each of its end vertices.
\end{lemma}

\begin{proof}
Let $\sigma_0 \sigma_1 \ldots \sigma_{n/2-1}$ be the disjoint cycle
decomposition of $s_m'$, where
\begin{equation*}
\sigma_i = ((m-i) \bmod n, (m+i+1) \bmod n) \quad (i=0,\dots,n/2-1).
\end{equation*}
Similar to the preceding proof, $\Gamma_n/s_m'$ has vertices
$\overline{\sigma_0},\dots,\overline{\sigma_{n/2-1}}$, and any two of them,
$\overline{\sigma_i}$ and $\overline{\sigma_j}$ (not necessarily distinct) are
adjacent in $\Gamma_n/s_m'$ if there exist adjacent vertices $v,w\in\Gamma_n$
such that $v\in\sigma_i$ and $w\in\sigma_j$; that is,
$v\equiv w\pm 1 \pmod{n}$, $v=(m-i)\bmod n$ or $v=(m+i+1) \bmod n$, and
$w=(m-j) \bmod n$ or $w=(m+j+1) \bmod n$.  The conjunction of these three
conditions is equivalent to $i=j\pm 1$ or $i=j=0$ or $i=j=n/2-1$. Therefore,
$\Gamma_n/s_m'$ consists of the path
$(\overline{\sigma_0},\overline{\sigma_1},\dots,\overline{\sigma_{n/2-1}})$
with loops at $\overline{\sigma_0}$ and $\overline{\sigma_{n/2-1}}$.
\end{proof}

In the following theorem, we use $\varphi(n)$ to denote Euler's totient
function, which gives the number of positive integers less than or equal to
$n$ that are coprime to~$n$, and $d\mid n$ to denote that $d$ is a positive
divisor of~$n$.  For $\mu=0$ the theorem specializes to Moreau's general
necklace polynomial \cite{Moreau:1872}, also known as the cycle index
polynomial of a cyclic group, while for $\mu=\lambda$ it specializes to the
orbital chromatic polynomial of a cycle of length~$n$.

Recall from Subsection~\ref{sec:conn-with-lucas} that $V_n(P,Q)$
denotes the Lucas sequence of the second kind.

\begin{theorem}
\label{thm:subgroup}
For every $n,\lambda\in\mathbb{N}$ and each $\mu = 0,\dots,\lambda$, we have
\begin{equation}
\label{eq:thmeq1}
OP_{\Gamma_n,\RG{n}}(\lambda,\mu) \, = \,\frac{1}{n} \sum_{d\mid n}
\varphi\!\left(\frac{n}{d}\right)
V_d(\lambda-1,\mu-\lambda)
\,-\,\frac{\mu-1}{2}\left( 1 - (-1)^n \right) .
\end{equation}
\end{theorem}

\begin{proof}
By Theorem~\ref{thm:orbital} and Eq.~\eqref{eq:4}, we have
\begin{equation}
\label{eq:subgroup}
OP_{\Gamma_n,\RG{n}}(\lambda,\mu) \,=\, \frac{1}{n} \sum_{m=0}^{n-1}
P_{\Gamma_n/r_m}(\lambda,\mu) .
\end{equation}
By Lemma \ref{lem:rotate},
$\Gamma_n/r_m = \Gamma_{\gcd(m,n)}$ if $\gcd(n.m)\neq 2$, and
$\Gamma_n/r_m = \Pi_1$ (a path of length~1)
if $\gcd(m,n)=2$. Since $P_{\Pi_1}(\lambda,\mu) = P_{\Gamma_2}(\lambda,\mu)$, we
conclude that 
\begin{equation*}
OP_{\Gamma_n,\RG{n}}(\lambda,\mu) \,=\,
\frac{1}{n} \sum_{m=0}^{n-1}
P_{\Gamma_{\gcd(m,n)}}(\lambda,\mu).
\end{equation*}
Using the identity $\gcd(m,n) = \gcd(n,n-m)$ and rearranging terms we obtain
\begin{equation*}
OP_{\Gamma_n,\RG{n}}(\lambda,\mu) \,=\, \frac{1}{n} \sum_{m=0}^{n-1}
P_{\Gamma_{\gcd(n,n-m)}}(\lambda,\mu) \,=\, \frac{1}{n} \sum_{m=1}^{n} P_{\Gamma_{\gcd(m,n)}}(\lambda,\mu)\, .
\end{equation*}
With $\varphi_d(n)=\#\{m\in\{1,\dots,n\}\mid \gcd(m,n)=d\}$ ($1\le d\le n$) we have
\begin{align*}
OP_{\Gamma_n,\RG{n}}(\lambda,\mu)  & =  \frac{1}{n} \sum_{d\mid n} \varphi_d(n)
P_{\Gamma_d}(\lambda,\mu) \,=\, \frac{1}{n} \sum_{d\mid n}
                                     \varphi\left(\frac{n}{d}\right) P_{\Gamma_d}(\lambda,\mu),
\end{align*}
and hence, by Eq.~(\ref{eq:51}),
\begin{align}
OP_{\Gamma_n,\RG{n}}(\lambda,\mu) & =  \frac{1}{n} \sum_{d\mid n}
                                    \varphi\!\left(\frac{n}{d}\right)
 \left( V_d(\lambda-1,\mu-\lambda) + (-1)^d(\mu-1) \!\right) \notag \\
  & =  \frac{1}{n} \sum_{d\mid n} \varphi\!\left(\frac{n}{d}\right)
    V_d(\lambda-1,\mu-\lambda)
    - \frac{\mu-1}{n} \sum_{d\mid n} (-1)^{d-1}
    \varphi\left(\frac{n}{d}\right).
    \label{eq:27}
\end{align}
As a consequence of Euler's divisor-sum identity, the second sum in
Eq.~\eqref{eq:27} equals $n$ if $n$ is odd, and $0$ if $n$ is even; formally,
\begin{equation}
\label{eq:29}
\sum_{d\mid n} (-1)^{d-1} \varphi\left(\frac{n}{d}\right) = \frac{n}{2}\left(
1 - (-1)^n \right).
\end{equation}
Eqs.~(\ref{eq:thmeq1}) now follows from Eq.~\eqref{eq:27} and
Eq.~\eqref{eq:29}.
\end{proof}

For $n=1,\dots,6$, the orbital bivariate chromatic polynomials given by
Theorem~\ref{thm:subgroup} are shown in Table~\vref{tab:orbital_cycle_Rot}.

\renewcommand{\arraystretch}{1.2}
\begin{table}[!t]
\footnotesize
\begin{center}
\begin{tabular}{rc}
  $n$ & $OP_{\Gamma_n,\RG{n}}(\lambda,\mu)$ \\[1pt] \hline
1 & ${\lambda} - {\mu}$ \\
2 & $\frac{1}{2} \, {\lambda}^{2} + \frac{1}{2} \, {\lambda} - {\mu}$ \\
3 & $\frac{1}{3} \, {\left({\lambda}^{2} - 3 \, {\mu} + 2\right)} {\lambda}$ \\
4 & $\frac{1}{4} \, {\lambda}^{4} - {\lambda}^{2} {\mu} + \frac{1}{4} \,
    {\lambda}^{2} + {\lambda} {\mu} + \frac{1}{2} \, {\mu}^{2} + \frac{1}{2}
    \, {\lambda} - \frac{3}{2} \, {\mu}$ \\
5 & $\frac{1}{5} \, {\lambda}^{5} - {\lambda}^{3} {\mu} + {\lambda}^{2} {\mu}
    + {\lambda} {\mu}^{2} - {\lambda} {\mu} - {\mu}^{2} + \frac{4}{5} \,
    {\lambda}$ \\
6 & $\frac{1}{6} \, {\lambda}^{6} - {\lambda}^{4} {\mu} + {\lambda}^{3} {\mu}
    + \frac{3}{2} \, {\lambda}^{2} {\mu}^{2} + \frac{1}{6} \, {\lambda}^{3} -
    {\lambda}^{2} {\mu} - 2 \, {\lambda} {\mu}^{2} - \frac{1}{3} \, {\mu}^{3}
    + \frac{1}{3} \, {\lambda}^{2} + \frac{1}{2} \, {\lambda} {\mu} +
    \frac{3}{2} \, {\mu}^{2} + \frac{1}{3} \, {\lambda} - \frac{7}{6} \,
    {\mu}$ 
\\[1pt] \hline
\end{tabular}
\caption{Orbital bivariate chromatic polynomials of $\Gamma_n$ with respect to
  $\RG{n}$}
\label{tab:orbital_cycle_Rot}
\end{center}
\end{table}

The following theorem expresses the orbital bivariate chromatic polynomial of
$\Gamma_n$ with respect to $\Aut(\Gamma_n)$ in terms of its orbital bivariate
chromatic polynomial with respect to $\RG{n}$. The bivariate chromatic
polynomial of a path, which appears in this expression, is given by
Eq.~\eqref{eq:6} or, equivalently, by Eq.~\eqref{eq:61}.

\begin{theorem}
\label{thm:fullgroup}
For every $n,\lambda\in\mathbb{N}$, with $n\ge 3$, and each $\mu =
0,\dots,\lambda$, we have
\begin{equation}
\label{eq:thmeq3}
OP_{\Gamma_n,\AG{n}}(\lambda,\mu) \,=\, \frac{1}{2} \,
OP_{\Gamma_n,\RG{n}}(\lambda,\mu) + \frac{\lambda-\mu}{2} P_{\Pi_{(n-3)/2}}(\lambda,\mu),
\end{equation}
if $n$ is odd, and
\begin{equation}
\label{eq:thmeq4}
OP_{\Gamma_n,\AG{n}}(\lambda,\mu)  = \frac{1}{2}
OP_{\Gamma_n,\RG{n}}(\lambda,\mu) + \frac{1}{4} P_{\Pi_{n/2}}(\lambda,\mu) +
\frac{(\lambda-\mu)^2}{4} P_{\Pi_{n/2-3}}(\lambda,\mu),                         
\end{equation}
if $n$ is even, where $\Pi_{-1}$ is interpreted as the empty graph.
\end{theorem}

\begin{proof}
For odd $n\ge 3$, by Theorem~\ref{thm:orbital} and Eq.~(\ref{ag1}),
\begin{equation*}
OP_{\Gamma_n,\AG{n}}(\lambda,\mu) \,=\, \frac{1}{2n} \left(
\sum_{m=0}^{n-1} P_{\Gamma_n/r_m}(\lambda,\mu) + \sum_{m=0}^{n-1}
P_{\Gamma_n/s_m}(\lambda,\mu) \right) .
\end{equation*}
By Eq.~(\ref{eq:subgroup}), the first sum agrees with
$n\, OP_{\Gamma_n,\RG{n}}(\lambda,\mu)$.  By Lemma
\ref{lem:reflect1}, $\Gamma_n/s_m$ is a path of length $(n-1)/2$ with a loop
attached to one of its end vertices.  For the color of this end vertex there
$\lambda-\mu$ choices, while for the remaining vertices there are
$P_{\Pi_{(n-3)/2}}(\lambda,\mu)$ choices.  Thus, in the second sum,
$P_{\Gamma_n/s_m}(\lambda,\mu) = (\lambda-\mu)P_{\Pi_{(n-3)/2}}(\lambda,\mu)$.
This proves Eq.~(\ref{eq:thmeq3}).
\par
For even $n\ge 4$, by Theorem~\ref{thm:orbital} and Eq.~(\ref{ag2}),
\begin{equation*}
OP_{\Gamma_n,\AG{n}}(\lambda,\mu) \,=\, \frac{1}{2n} \left( 
\sum_{m=0}^{n-1} P_{\Gamma_n/r_m}(\lambda,\mu) + \sum_{m=0}^{n/2-1}
P_{\Gamma_n/s_m}(\lambda,\mu) + \sum_{m=0}^{n/2-1}
P_{\Gamma_n/s_m'}(\lambda,\mu) \right) .
\end{equation*}
Again, the first sum agrees with $n\, OP_{\Gamma_n,\RG{n}}(\lambda,\mu)$.  By
Lemma \ref{lem:reflect1}, $\Gamma_n/s_m$ is a path of length $n/2$,
whence $P_{\Gamma_n/s_m}(\lambda,\mu) = P_{\Pi_{n/2}}(\lambda,\mu)$.  By
Lemma \ref{lem:reflect2}, $\Gamma_n/s_m'$ is a path of length $n/2-1$
with a loop attached to both end vertices. Similar to the odd case,
$P_{\Gamma_n/s_m'}(\lambda,\mu) = (\lambda-\mu)^2P_{\Pi_{n/2-3}}(\lambda,\mu)$
if $n>4$, and $P_{\Gamma_n/s_m'}(\lambda,\mu) = (\lambda-\mu)^2$ if
$n=4$, which proves Eq.~(\ref{eq:thmeq4}).
\end{proof}

For $n=1,\dots,6$, the orbital bivariate chromatic polynomials provided by
Theorem~\ref{thm:fullgroup} are shown in Table~\vref{tab:orbital_cycle_Aut}.

\renewcommand{\arraystretch}{1.2}
\begin{table}[!t]
\footnotesize
\begin{center}
\begin{tabular}{rc}
  $n$ & $OP_{\Gamma_n,\AG{n}}(\lambda,\mu)$ \\[1pt] \hline
1 & ${\lambda} - {\mu}$ \\
2 & $\frac{1}{2} \, {\lambda}^{2} + \frac{1}{2} \, {\lambda} - {\mu}$ \\
3 & $\frac{1}{6} \, {\left({\lambda}^{2} + 3 \, {\lambda} - 6 \, {\mu} +
    2\right)} {\lambda}$ \\
4 & $\frac{1}{8} \, {\left({\lambda}^{2} + {\lambda} - 2 \, {\mu} + 2\right)}
    {\left({\lambda}^{2} + {\lambda} - 2 \, {\mu}\right)}$ \\
5 & $\frac{1}{10} \, {\left({\lambda}^{4} - 5 \, {\lambda}^{2} {\mu} + 5 \,
    {\lambda}^{2} + 5 \, {\mu}^{2} - 10 \, {\mu} + 4\right)} {\lambda}$ \\
6 & $\frac{1}{12} \, {\lambda}^{6} - \frac{1}{2} \, {\lambda}^{4} {\mu} +
    \frac{1}{4} \, {\lambda}^{4} + \frac{1}{2} \, {\lambda}^{3} {\mu} +
    \frac{3}{4} \, {\lambda}^{2} {\mu}^{2} + \frac{1}{3} \, {\lambda}^{3} -
    \frac{7}{4} \, {\lambda}^{2} {\mu} - \frac{3}{4} \, {\lambda} {\mu}^{2} -
    \frac{1}{6} \, {\mu}^{3} + \frac{1}{6} \, {\lambda}^{2} + \frac{3}{4} \,
    {\lambda} {\mu} + {\mu}^{2} + \frac{1}{6} \, {\lambda} - \frac{5}{6} \,
    {\mu}$ 
\\[1pt] \hline
\end{tabular}
\caption{Orbital bivariate chromatic polynomials of $\Gamma_n$ with respect to
  $\AG{n}$}
\label{tab:orbital_cycle_Aut}
\end{center}
\end{table}

\subsection{Wheels}

The wheel $\Wheel_n$ on $n\ge 3$ vertices is the join of $\Kappa_1$ with
$\Gamma_{n-1}$. In the following, we consider the case $n\ge 5$, where
$\Aut(\Wheel_n) \cong \Aut(\Gamma_{n-1})$.

\begin{theorem}
Let $n\in\mathbb{N}$, $n\ge 5$, and $G$ be a subgroup of $\Aut(\Wheel_n)$.
Then, for every $\lambda\in\mathbb{N}$ and $\mu=0,\dots,\lambda$ we have
\[ OP_{\Wheel_n,G}(\lambda,\mu)
= \mu \,{OP}_{\Gamma_{n-1},G}(\lambda-1,\mu-1) +
(\lambda-\mu)\,{OP}_{\Gamma_{n-1},G}(\lambda,\mu) . \]
\end{theorem}

\begin{proof}
Since each $g\in G$ fixes the central vertex, the non-equivalent $\mu$-proper
$\lambda$-colorings can be counted by a case distinction: If the central
vertex receives a color from $\{1,\dots,\mu\}$, the remaining $n-1$ vertices
(forming the cycle $\Gamma_{n-1}$) must be colored using the remaining
$\lambda-1$ colors, with $\mu-1$ colors chosen from $\{1,\dots,\mu\}$.  This
yields $\mu\,OP_{\Gamma_{n-1},G}(\lambda-1,\mu-1)$ non-equivalent colorings.
If the central vertex receives a color from $\{\mu+1,\dots,\lambda\}$, this
imposes no restriction on the colors of the vertices of the cycle, which
yields $(\lambda-\mu)\,OP_{\Gamma_{n-1},G}(\lambda,\mu)$ non-equivalent
colorings in this case. Summing these two cases gives the result.
\end{proof}

As an example, we note the orbital bivariate chromatic polynomial of
$\Wheel_5$ with respect to its automorphism group:
\begin{multline*}
OP_{\Wheel_5, \Aut(\Wheel_5)} (\lambda,\mu) = \frac{1}{8} \, {\lambda}^{5} + \frac{1}{4} \, {\lambda}^{4} - {\lambda}^{3} {\mu} + \frac{3}{8} \, {\lambda}^{3} + \frac{3}{2} \, {\lambda} {\mu}^{2} + \frac{1}{4} \, {\lambda}^{2} - \frac{3}{2} \, {\lambda} {\mu} - {\mu}^{2} + {\mu}
\end{multline*}

\section{A final treat}
\label{sec:final-treat}

\subsection{Generalization of Moreau's necklace polynomial}

Theorem~\ref{thm:subgroup} provides a bivariate generalization of Moreau's
necklace polynomial. Moreau's polynomial counts equivalence classes of
periodic sequences (necklaces) under rotation where any color combination of
adjacent beads from the set of colors $\{1,\dots,\lambda\}$ is permitted,
which corresponds to the case $\mu=0$ in Theorem~\ref{thm:subgroup}.  Clearly,
by Eq.~\eqref{eq:25} from Subsection~\ref{sec:conn-with-lucas},
$V_n(\lambda-1,-\lambda) = \lambda^n + (-1)^n$ for every
$n,\lambda\in\mathbb{N}$, whence by Eqs.~\eqref{eq:thmeq1} and~\eqref{eq:29},
\begin{align*}
OP_{\Gamma_n,\RG{n}}(\lambda,0) & = \frac{1}{n} \sum_{d\mid n}
\varphi\!\left(\frac{n}{d}\right) \left( \lambda^d + (-1)^d \right) +
                                    \frac{1}{2} \left( 1-(-1)^n\right) 
  \\
  & = \frac{1}{n} \sum_{d\mid n}
\varphi\!\left(\frac{n}{d}\right) \lambda^d  - \frac{1}{n} \sum_{d\mid n}
    (-1)^{d-1} \varphi\!\left(\frac{n}{d}\right) + \frac{1}{2} \left(
    1-(-1)^n\right) \\
  & = \frac{1}{n} \sum_{d\mid n}
    \varphi\!\left(\frac{n}{d}\right) \lambda^d \, ,
\end{align*}
which reveals Moreau's original necklace polynomial,

The case $\mu>0$ is more interesting, as it imposes restrictions on the colors
of adjacent beads. Adjacent beads may be assigned the same color only if that
color is chosen from a distinguished subset of \(\lambda-\mu\) colors. In the
special case \(\mu=\lambda\), this means that adjacent beads must receive
distinct colors.

\subsection{Generalization and proof of Fermat’s Little Theorem}

As a consequence of Theorem~10, we obtain a generalization of Fermat's little
theorem for the Lucas sequences $V_n(P,Q)$, as defined in
Subsection~\ref{sec:conn-with-lucas}. Restricted versions were obtained by
Lucas~\cite{Lucas:1878}, assuming that $P$ and $Q$ are coprime and 
$n$ divides $P^2 - 4Q$, and by Carmichael~\cite{Carmichael:1913}, assuming only
that $P$ and $Q$ are coprime and, unlike in our proof, relying on Fermat's
little theorem itself.  Notably, our generalization implies Fermat's little
theorem.

\begin{theorem}
\label{cor:vpq}
For every prime number $n$ and every $P,Q\in\mathbb{Z}$ we have
\begin{gather*}
V_n(P,Q) \equiv P \pmod{n} .
\end{gather*}
\end{theorem}

\begin{proof}
Putting $\lambda = (P\bmod n)+n+1$ and $\mu=(P \bmod n) + (Q\bmod n)+1$, we
have $1\le\mu\le \lambda-1$, $P\equiv \lambda-1 \pmod{n}$, and
$Q\equiv \mu-\lambda \pmod{n}$. Thus, by Theorem~\ref{thm:subgroup},
\begin{gather*}
\sum_{d\mid n} \varphi\left( \frac{n}{d} \right) V_d(P,Q) \equiv  0 \pmod{n}
.
\end{gather*}
Since $n$ is a prime number, this congruence reduces to
\[ V_n(P,Q) \equiv V_1(P,Q) = P \pmod{n}, \]
which proves the result.
\end{proof}

Theorem~\ref{cor:vpq} implies Fermat's Little Theorem.

\begin{corollary}[Fermat's Little Theorem]
For every prime number $n$ and every $m\in\mathbb{Z}$,
\begin{gather*}
m^n \equiv m \pmod{n}.
\end{gather*}
\end{corollary}

\begin{proof}
By Eq.~\eqref{eq:25} and Theorem~\ref{cor:vpq} we have
$m^n = V_n(m,0) \equiv m \pmod{n}$, proving the statement. Alternatively, for
$n>2$, we have $m^n-1 = V_n(m-1,-m) \equiv m-1\pmod{n}$.
\end{proof}

\printbibliography

\end{document}